# A Two-Stage Metaheuristic Algorithm

# for the Dynamic Vehicle Routing Problem

# in Industry 4.0 approach

Maryam Abdirad[1], Krishna Krishnan[1], Deepak Gupta[1]

## Abstract

Industry 4.0 is a concept that assists companies in developing a modern supply chain (MSC) system when they are faced with a dynamic process. Because Industry 4.0 focuses on mobility and real-time integration, it is a good framework for a dynamic vehicle routing problem (DVRP). This research works on DVRP. The aim of this research is to minimize transportation cost without exceeding the capacity constraint of each vehicle while serving customer demands from a common depot. Meanwhile, new orders arrive at a specific time into the system while the vehicles are executing the delivery of existing orders. This paper presents a two-stage hybrid algorithm for solving the DVRP. In the first stage, construction algorithms are applied to develop the initial route. In the second stage, improvement algorithms are applied. Experimental results were designed for different sizes of problems. Analysis results show the effectiveness of the proposed algorithm.

**Keywords:** Dynamic Vehicle Routing Problem, Industry 4.0, Two-Stage algorithm, Heuristic algorithms

Maryam Abdirad
Ph.D. Candidate at Wichita State University
mxabdirad@wichita.edu
[1] Department of Industrial, Systems, and Manufacturing Engineering, Wichita State University, Wichita, KS 67260, USA

# 1.    Introduction

Achieving a modern and agile supply chain (SC) that is efficient, automated, flexible, and transparent is the goal of most companies. Moreover, a modern supply chain (MSC) can work in a dynamic system and can handle high volumes of data, enabling efficient cooperation among all elements of the SC including suppliers, manufacturers, and customers [1][2]. An outstanding example of MSC is Amazon company which is a company, that provides a fast response to customers in preparing, shipping, and delivering the product to the customers. Amazon uses the combination of supply chain network with Industry 4.0 to make this company unique among other companies. Industry 4.0 provides a framework that can guide the move from a traditional SC to an MSC. This strategic approach focuses on automation, digitalization, interconnection (e.g., via the Internet of Things [IoT]), information transparency, and decentralized decisions (e.g., autonomous cyber-physical systems) in companies. Industry 4.0 focuses on mobility and real-time integration, and hence it can provide a good framework for the SC problem [2][3].

One of the well-known supply chain problems is the vehicle routing problem (VRP), which looks for an optimal set of routes to deliver demands to demand points. Different variations of the VRP take into account several features of this problem [4], such as Capacitated VRP (CVRP), Multi-Depot VRP (MDVRP) and VRP with time windows (VRPTW). The one that has most recently received considerable attention is the dynamic vehicle routing problem (DVRP). The DVRP is a real-world problem with high complexity and intractable nature which has to be solved for dynamic supply chain systems and is also referred to as an online or real-time vehicle routing problem [5]. The transportation and logistics problems are optimized using a static model, but with the increase in traffic and demand along with the demand for flexibility by customers, there is an increase in computational and communication needs for solving the DVRP in dynamic conditions. As Industry 4.0 can handle a dynamic system, it can also be a good framework for the DVRP.

The main goal of this study to introduce the dynamic vehicle routing problem with a single depot and develop a two-stage algorithm to solve it. The DVRP has dynamic demands from customers at different locations that arrive in the system at different times. Each new demand and customer obviously affects the solution because they change both the problem and the solution the instant they arrive in the system. The challenge of this research and the objective here is the construction of routes from a depot with minimum distances to the destination. In this paper, we formulate the dynamic vehicle routing problem (DVRP) as an integer program model to make the routing decision.

The remaining parts of this paper are organized as follows. Section 2 provides a brief literature review dedicated to Industry 4.0 and its role in the supply chain. This is followed by a brief review of Industry 4.0 and the vehicle routing problem, with an emphasis on DVRPs. The problem is described in Section 3. In section 4, the solution approach to this problem is explained. Three different scenarios are presented in the experimental results section, and the last section concludes with a summary and an outlook for future work.

## 2.    Literature Review

### 2.1.    *Industry 4.0*

The concept of Industry 4.0 was presented in 2011 by Henning Kagermann (former top manager of the German company SAP) [6]. Industry 4.0 is alternatively known as the "Fourth Industrial Revolution," "Smart Manufacturing," "Industrial Internet," or "Integrated Industry"[7]. This concept is becoming increasingly more popular and has been receiving attention all over the world [8][9]. However, industry experts have not identified a precise definition of Industry 4.0.

According to Lopes de Sousa Jabbour et al. [10], "the core feature of Industry 4.0 is connectivity between machines, orders, employees, suppliers, and customers for the Internet of Things and electronic devices. Consequentially, firms can produce products using decentralized decisions and autonomous systems".

The first industrial revolution began with the development of water power and steam power and the mechanization of the production system in 1784. The second industrial revolution changed the production system to a mass production system and advanced assembly lines by the use of electricity in the 1870s. The third industrial revolution was a big revelation to automate some of the production processes by using computers in 1970. The fourth industrial revolution leads all integrations of a system to digitalization by using IoT and cyber-physical systems (CPSs), termed "Industry 4.0" [11][12]. Armengaud et al. defined Industry 4.0 as [13] "the comprehensive introduction of information and communication technology (ICT) as well as their connection to an internet of things, services and data, which enables a real-time production. Industry 4.0 means a higher degree of digitalization for products, value creation chain, and business models". Figure 1 demonstrates the four industrial revolutions.

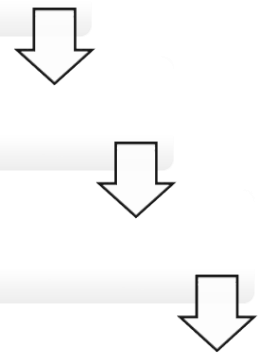

**First Industrial Revolution**
• Mechanization, Water and Steam Power Engine (1784)

**Second Industrial Revolution**
• Mass Production, Assembly Line using Electriacal Energy (1870)

**Third Industrial Revolution**
• Use of PLC and IT systems for Automation (1970)

**Fourth Industrial Revolution**
• Use of IoT and Cyber Physical System (Today)

Figure 1. The Four Industrial Revolutions[14].

The prime focus of Industry 4.0 is to have a smart network based on digitalization and automatization where machines and products interact with each other without human involvement [15][16]. The outcome of Industry 4.0 is the development of smart factory systems that include smart machines, smart devices, smart manufacturing processes, smart engineering, smart logistics, smart suppliers, smart products, etc. [17][18][19]. Industry 4.0 promotes the use of CPS, the Internet of Things, the Internet of Services (IoS), robotics, big data, and cloud manufacturing. Thus, devices, machines, production modules and products, etc. are applied to various fields, such as the supply chain, manufacturing, and management, especially in response to real-time situations [20][21][22]. The interested readers referred to read about CPS [23][24][25][26].

Industry 4.0 is expected to have a significant impact on supply chains, business models, and processes related to achieving a modern supply chain. Researchers use different names for Industry 4.0 in the supply chain, such as digital supply network (DSN), Internet of Things (IoT), electronic supply chain (E-supply chain), Supply Chain 4.0, E-logistic, or Logistic 4.0. As explained previously, Industry 4.0 increases digitalization and automation in manufacturing and creates a digital process to facilitate interaction among all parts of the company. By implementing Industry 4.0 in the supply chain system, four main supply chain elements—integration, operations, purchasing, and distribution—are affected, which can also increase the productivity of companies [27][3]. The main benefits of Industry 4.0 in the supply chain are to reduce lead time for the delivery of products to customers, to reduce the time needed to respond to an unforeseen event, and to prompt a significant increase in the quality of decision-making [28]. Industry 4.0 can help companies afford complicated and dynamic processes in their supply chain and to handle large-scale production and integration of customers [9]. Industry 4.0 can bring positive benefits in current sales and operations planning and also in the logistics process [29]. After implementing Industry 4.0, real-time information can be shared across this digitalized process to drive useful decisions.

In Industry 4.0 in SC, there is communication between systems, including the supply chain management (SCM) control tower, depot, and the vehicles (Figure 2). Technological advancements, such as mobile devices, enable direct communication between them. So, a driver can dynamically change his/her plan while executing the route. Also, the emergence of the global positioning system (GPS) allows the SCM control tower to know the current position of a driver and communicate in a timely manner the next customer to visit on the route [5]. Additionally, using radio frequency identification (RFID) chips and sensors in packages can help to facilitate this type of communication.

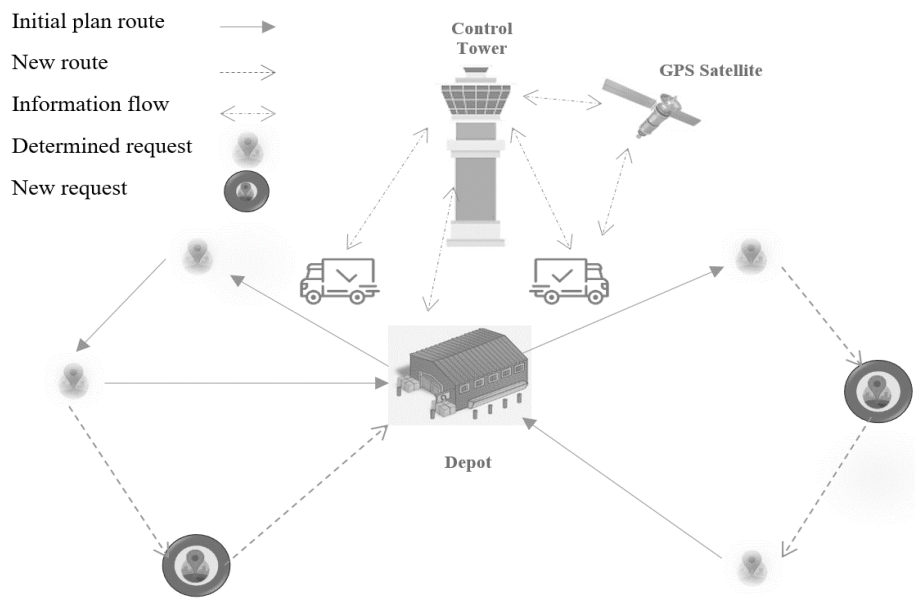

Figure 2. Communication network for real-time measurements of DVRP [30].

## 2.2. *Vehicle Routing Problem*

The vehicle routing problem is not a new subject and has been studied for over five decades. The formulation for VRP was first introduced by Dantzig and Ramser as part of the Traveling Salesman

Problem (TSP)[31][32]. VRP can be described as finding an optimal expected cost for delivery from one or many depots to many customers who are geographically distributed. It assumes that the travel distance between customers is the Euclidean distance between their coordinate pairs. The cost is considered to be proportional to the distance traveled. Each customer must be assigned to exactly one of the K vehicle routes, and the total demands assigned to each vehicle must not exceed the vehicle capacity. VRP is an NP-hard problem (non-deterministic polynomial time) and is considered a combinatorial optimization and integer programming problem.

VRP is categorized into four categories: static and deterministic (SD), static and stochastic (SS), dynamic and deterministic (DD), and dynamic and stochastic (DS) [33]. The focus of this research is on the DD vehicle routing problem. Wilson and Colvin presented the first paper about DVRP [34]. The DD problem is also referred to as an online or real-time problem. All data related to the routing process are not known before planning, and they can change during the planning horizon. More explanation about the DD and DS categories are provided below [31][33]:

- Dynamic and deterministic: In this category, all inputs are unknown and revealed dynamically during the design or execution of the routing plans. In this setting, the information is stochastic and hence future information is not known. For example, the location of a customer may be unknown until that customer request is received.

- Dynamic and stochastic: In the dynamic and stochastic category, parts or all inputs are unknown and revealed dynamically during the execution of the routes. However, in contrast to the DD problem, in addition to efficiently handling dynamic events, stochastic knowledge about the revealed data is also available.

Figure 3 illustrates the difference between DVRP and the classic problem. In the beginning, when time = 0, there is no customer in the system. Customer A arrives into the system at time = 5.

Customer A is satisfied between time 5 and time 10. Customer B arrives at time = 10. Customer B is satisfied between time 10 and 15. Customer C arrives at time t = 15, but there are not enough products in the vehicle to deliver to customer C. The vehicle must return to the depot to refill (time = 20). Thus, customer C is serviced at time t = 25, since customers A and B are already serviced [31].

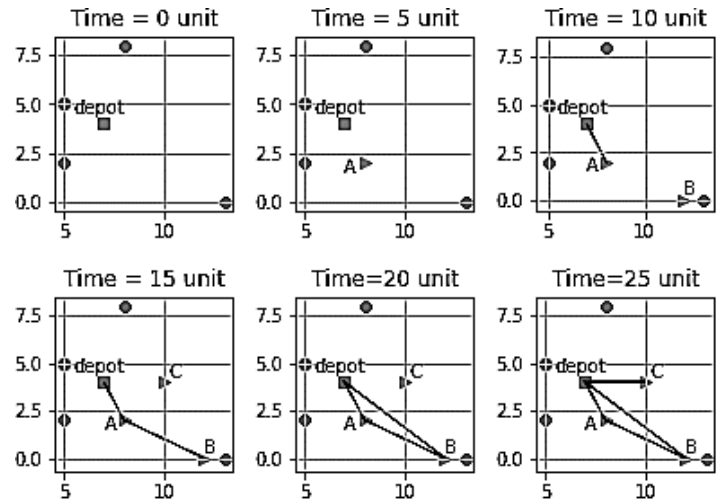

Figure 3. Example of dynamic vehicle routing.

Through this introduction, dynamic events can be categorized into three different classes [35]:

- Dynamically incoming or outgoing entities, e.g., orders or vehicles.

- Dynamically occurring unexpected events.

- New or changing information transmitted dynamically.

The first published paper on the dynamic vehicle routing problem involved the dial-a-ride-problem (DARP) where customers appear dynamically in the system [21]. After that, in 1980, the concept of an immediate request was published by Psaraftis et al., where the current route has to be changed by coming a new customer request in order to respond to [36]. Following this research, a different approach used to find the best solution.

Exact algorithms are those algorithms or optimization models that can achieve an optimal solution [33]. The first published paper involving exact algorithms to the VRP was based on branch-and-bound, set partitioning, and branch-and-cut [24]. One approach of exact algorithms is the branch-and-cut algorithm. Desaulniers proposed a new branch-and-price-and cut algorithm for the exact solution of split-delivery vehicle routing problem with time windows [37]. Also, Toth and Vigo proved that the branch-and-cut-and-price algorithm is a good exact method for the VRP [38]. A branch-and-bound algorithm is an exact algorithm for the vehicle routing problem based on a spanning tree and shortest-path relaxations. Ropke and Cordeau introduced a branch-and-cut algorithm for pickup and delivery problems with time windows. In this algorithm, lower bounds are computed by solving through column generation the linear programming relaxation of a set partitioning formulation [39].

Most of the time, real-world vehicle routing problems are large. It is difficult to solve the VRP with exact solutions, especially in a reasonable period of time. Thus, different types of heuristic methodologies can be a good approach. The Clarke and Wright (CW) Savings algorithm was applied for the first time in the VRP [27]. Dror and Trudeau modified savings algorithm to illustrate the effects of route failure on the expected cost of a route [40]. The sweep (SW) algorithm is a method for clustering customers into groups so that customers in the same group are geographically close together and can be served by the same vehicle [23]. Nurcahyo, Alias, and Shamsuddin applied the SW algorithm in solving a VRP for public transportation [41]. The application of particle swarm optimization (PSO) and the performance result in the DVRP was shown by Khouadjia et al. [42]. Kergosien et al. applied Tabu Search (TS) to solve the transportation of patients in a hospital. Some demands are known, and others enter the system dynamically [43]. Gendreau et al. used the neighborhood search (NS) algorithm to solve the dynamic pickup and delivery problem with new customers [44]. Elhassania, Jaouad and Ahmed reported using a genetic algorithm (GA) for the DVRP [45]. Also, Hanshar and Ombuki-Berman used a GA for providing a solution for DVRP [46].

Some papers proposed new algorithms for solving the DVRP, for instance, ant colony optimization (ACO) algorithms [47], enhanced ant colony optimization (ACO) [48], optimization utilizing monarch butterfly optimization [49], and firefly algorithm [50]. Elhassania et al. proposed a hybridization obtained by combining an ACO algorithm with a large neighborhood search (LNS) algorithm to solve DVRP with static demands [51]. In another paper by Novoa and Storer, an approximate dynamic programming approach was applied for the vehicle routing problem with stochastic demands [52].

## 3.     Problem

In this research, 'm' homogeneous vehicles with fixed equal capacity ($q_i$) and i = 1, . . ., m, depart from a depot to deliver products to 'n' customers at demand points. Each customer has a known demand $d_i$ (i = 1, . . ., n). It is assumed that the customers' demands are less than the maximum capacity of the vehicles. Meanwhile, new customers with known demand emerge dynamically over time. The distance of the route, calculated by assuming Euclidian distance, is associated with every edge in the total route. The solution may use all of the 'm' vehicles or a subset of the vehicles based on the available demands. Other constraints are given as follows:

- Each vehicle starts and ends its route at the depot.

- All customer demand ($d_i$ where i = 1, . . ., n) should be accepted.

- All customer demand must be satisfied.

- Each customer is assigned to be served by only one vehicle.

- The sum of the demands in each vehicle route does not exceed the vehicle's capacity.

- The unit of distance is identical to the unit of cost.

## 4.    Solution Approach

In VRPs, an exact algorithm may provide a solution in a reasonable time period. However, when modeling and solving DVRPs, it is often impossible to solve the DVRP with an exact algorithm in a reasonable time. Therefore, a metaheuristic approach is a good solution for this kind of problem. One of the methodologies that can be effectively used for the DVRP problem is a construction-route first, improvement-route second approach. To create the initial plan to deliver the product from the depot to the known customers, three different algorithms are executed. The two-stage algorithm (Figure 4) includes construction algorithms and improvement algorithms that are explained below:

*Construction Algorithms:* The goal of this stage is to construct a route for each generated cluster separately and to find the best overall routing objective with a feasible solution. Heuristic algorithms are used to get a solution in a reasonable time. In this research, Path Cheapest Arc, Savings, and Global Cheapest Arc are applied for the construction phase of route identification. In this research, the goal is to minimize the transportation cost.

*Improvement Algorithms:* The route obtained from the first stage construction algorithms may not be optimal and can be improved. The sub-optimal solution from the first stage is then fed into the second stage of the algorithm. The second-level routing heuristic is used to improve the solution obtained from the first stage. Therefore, improvement algorithms may be used to further improve the solution. There are different algorithms that can be used for the solution improvement stage. In this research, Guided Local Search, Simulated Annealing, and Tabu Search are used. To get a fast result, there is a possibility to add run time constraint for second stage algorithms

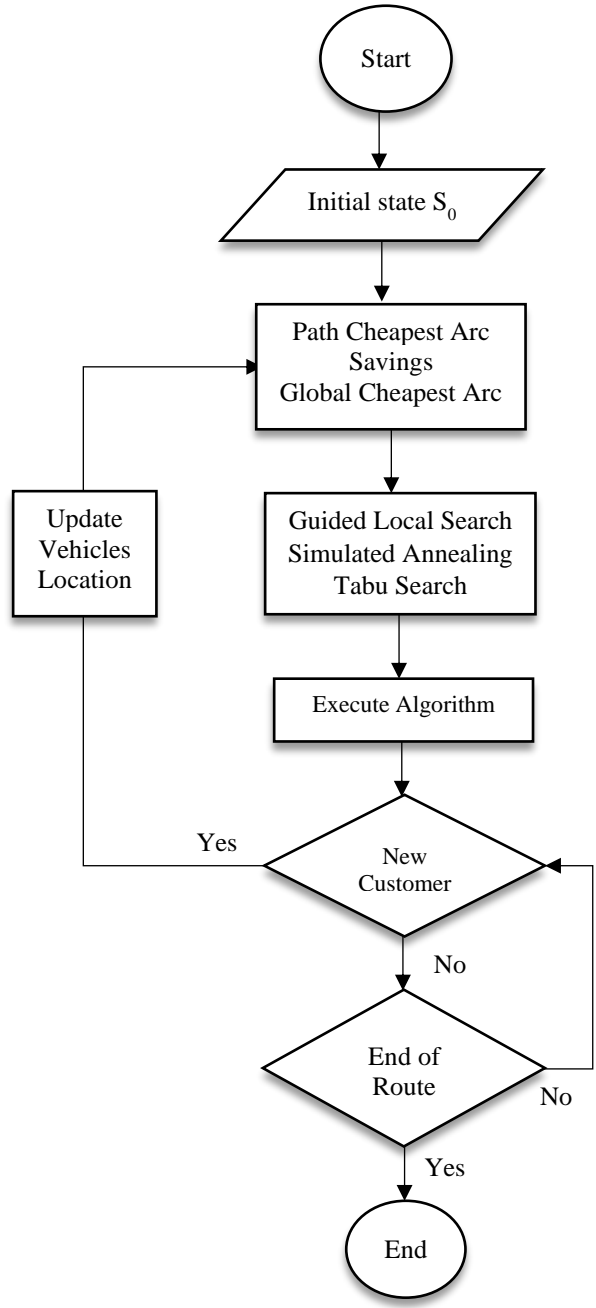

Figure 4. Details of applied methodology for DVRP.

At time $T = 0$, after executing the two-stage algorithm and the routes are determined, vehicles leave the depot to serve the customers. The initial routes are modified when new customers enter the system. If exchanging or inserting is done between and during the routes, it is called "inter-route improvement"[53]. Therefore the two-stage algorithm is executed each time a new customer appears until all customers are served [54].

## 5.    Experimental Results and Discussion

To verify the efficiency of the proposed method, several experiments are carried out. Then, experimental results for three different case studies in small, medium, and large sizes are shown. At the end, analysis of the results from case studies are discussed.

### *5.1.    Data Collection and Processing*

The proposed method is tested with nine different data sets including the different numbers of stat data and dynamic data. The algorithm was implemented using Python. Experiments are performed on a personal PC Intel® Core ™ i7- 4790S CPU @ 3.20 GHz, 3201 MHz, 4 Core(s), with 8 GB of RAM. To account for variability in the solution obtained, the two-stage algorithm is executed 10 times based on the calculation of the deviation in the solution at 95% confidence interval. The maximum standard deviation observed in all cases is 4.4.

### *5.2.    Experimental Results*

The proposed methodology was applied for all the dynamic delivery routing and the results were obtained. Table 1 gives the computational results for the proposed approach for the dynamic delivery routing instances. For each instance, the best and average solution and name of the algorithm that returned the best results were provided. To verify the effectiveness and efficiency of the proposed algorithm, the results generated with the Two-stage algorithm were compared with the first stage of the selected algorithm. The last column of Table 1 gives the performance comparison of the proposed approach for the dynamic delivery routing instances.

Table 1. Computational results for the proposed methodology for the dynamic delivery instances.

| Dataset | Stat Data | Dynamic Data | Average Cost ($) | Minimum Cost ($) | Two-Stage Algorithm | First stage Algorithms | Improvement (%) |
|---|---|---|---|---|---|---|---|
| Small size | 20 | 20 | 304.3 | 298.52 | Path Cheapest Arc and Tabu Search | 398.5 | 23.63 |
| | 20 | 50 | 532.76 | 523.79 | Path Cheapest Arc and Simulated Annealing | 616.04 | 13.51 |
| | 20 | 100 | 788.94 | 788.63 | Global Cheapest Arc and Tabu Search | 930.1 | 15.17 |
| Medium size | 50 | 20 | 542.61 | 533.27 | Savings and Tabu Search | 588.41 | 7.78 |
| | 50 | 50 | 660.71 | 649.32 | Savings and Tabu Search | 726.65 | 9.07 |
| | 50 | 100 | 865.39 | 849.94 | Global Cheapest Arc and Simulated Annealing | 921.55 | 6.09 |
| Large size | 100 | 20 | 990.20 | 983.57 | Global Cheapest Arc and Guided Local Search | 1080.91 | 8.39 |
| | 100 | 50 | 1123.03 | 1101.98 | Global Cheapest Arc and Guided Local Search | 1228.74 | 8.61 |
| | 100 | 100 | 1437 | 1431.37 | Global Cheapest Arc and Tabu Search | 1599.2 | 10.14 |

As evident from Table 1, the two-stage algorithm provides better results compared to the construction algorithm used in the first-stage. In summary, it was found that in small size instances with 20 initial static customers at time T=0 and 20, 50 and 100 dynamic customers who enter the system during the execution of the route, Path Cheapest Arc and Tabu Search, Path Cheapest Arc and Simulated Annealing, and Global Cheapest Arc and Tabu Search have the lowest transportation cost respectively. In medium-size instances with 50 initial customers, 20 and 50 dynamic customers, the lowest transportation cost is obtained when using Savings and Tabu Search algorithm. When the dynamic customer number is increased to 100, Global Cheapest Arc and Simulated Annealing have the lowest transportation cost. In large-size instances, with 100 initial customers and 20, 50 and 100 dynamic customers, Global Cheapest Arc and Guided Local Search, Global Cheapest Arc and Guided Local Search, and Global Cheapest Arc and Tabu Search have the lowest transportation cost, respectively.

In developing these case studies, it has been found that some case study data sets resulted in infeasible solutions in the construction algorithm phase. Further research into identifying the best construction algorithms based on the demand data patterns must be investigated. The two-stage algorithm results in either reduced transportation costs or the same as the construction algorithm transportation costs.

DVRP is an NP-Hard problem, and hence it is often hard and even impossible to reach an exact solution for real-world problems. Furthermore, for NP-hard problems, heuristic methods may also not be able to find local optimal solutions within a specific run time. One of the advantages of the two-stage algorithm proposed in this research is that it makes this problem solvable. One disadvantage of this method is that the achievement results are considered as the local optimum, and may not be the global optimum. The second drawback of this approach is that the heuristic must be executed before updating the solution, which can increase delays for the vehicles, while computational power is unused during waiting times

It has been found that by increasing run time, there may be fluctuations in transportation costs. However, in most cases, longer run times results in moving local optimal solutions to the global optimal solutions. As an outcome, the results indicate that the developed heuristic performs well and provides a good result on DVRP.

6.     **Conclusion and Future Work**

The present study solved the DVRP with a depot and a set of initial and dynamic customers whose demand should be fully satisfied. The delivery of products to the customer should be ensured by a single-vehicle. New demands appear in the system over time which make this problem dynamic. The objective of this problem includes the minimization of the total travel cost.

This problem is important in both the research and industrial domains due to its many real-world applications.

In this paper, a two-stage algorithm was proposed for solving a DVRP. In the first stage of the proposed methodology, construction algorithms such as Savings algorithm, Path cheapest Arc algorithm, and Global Cheapest Arc are used to construct the initial route. In the second stage, some improvement algorithms such as Guided Local Search, Simulated Annealing, and Tabu Search algorithms are applied to improve the initial route. To verify the proposed methodology, the new solutions are compared to existing metaheuristics in the same category. The lowest transportation cost from all these combinations is selected as the best answer to this problem. The two-stage algorithm was tested on multiple case studies to evaluate its effectiveness. As explained in the results section, the proposed two-stage algorithm results in lower transportation costs.

It is necessary to emphasize that this methodology may lead to better solutions or even the best solutions. The algorithms are executed for specific run times. As DVRP is a real-world problem, the two-stage algorithms can be applied to newer and larger case studies to determine its effectiveness in solving large-size problems. Considering that the execution times are relatively small, it should be possible to execute large size data sets. The research could also be extended to include the concept of hard and soft time-windows. For future research, an interesting approach would be to use this algorithm under the assumption that customer demand should be delivered in a specific time window. Besides, the use of loading and unloading times as additional inputs to the case study can also be developed. As the last suggestion for further study, the methodology can also focus on fuzzy data with probabilistic demands.

# References

[1]  Y. Yin, K. E. Stecke, and D. Li, "The evolution of production systems from Industry 2.0 through Industry 4.0," *Int. J. Prod. Res.*, vol. 56, pp. 848–861, 2018.

[2]  J. Barata, P. Rupino Da Cunha, and J. Stal, "Mobile supply chain management in the Industry 4.0 era: An annotated bibliography and guide for future research," *J. Enterp. Inf. Manag.*, vol. 31, no. 1, pp. 173–192, 2017.

[3]  M. Abdirad and K. Krishnan, "Industry 4.0 in Logistics and Supply Chain Management: A Systematic Literature Review," *EMJ - Eng. Manag. J.*, 2020.

[4]  K. Braekers, K. Ramaekers, and I. Van Nieuwenhuyse, "The vehicle routing problem: State of the art classification and review," 2016.

[5]  R. Necula, M. Breaban, and M. Raschip, "Tackling dynamic vehicle routing problem with time windows by means of ant colony system," in *2017 IEEE Congress on Evolutionary Computation (CEC)*, 2017, pp. 2480–2487.

[6]  W. Paprocki, "How transport and logistics operators can implement the solutions of 'Industry 4.0,'" in *Sustainable Transport Development, Innovation and Technology*, 2016, pp. 185–196.

[7]  E. Hofmann and M. Rüsch, "Industry 4.0 and the current status as well as future prospects on logistics," *Comput. Ind.*, vol. 89, pp. 23–34, 2017.

[8]  Y. Liao, F. Deschamps, E. de Freitas Rocha Loures, and L. Felipe Pierin Ramos, "Past, present and future of Industry 4.0—A systematic literature review and research agenda proposal," *Int. J. Prod. Res.*, vol. 55, no. 12, pp. 3609–3629, 2017.

[9]  F. Rennung, C. T. Luminosu, and A. Draghici, "Service provision in the framework of Industry 4.0," *Procedia - Soc. Behav. Sci.*, vol. 221, pp. 372–377, Jun. 2016.

[10]  A. B. Lopes de Sousa Jabbour, C. J. C. Jabbour, M. Godinho Filho, and D. Roubaud, "Industry 4.0 and the circular economy: A proposed research agenda and original roadmap for sustainable operations," *Ann. Oper. Res.*, pp. 1–14, 2018.

[11]  L. Li, "China's manufacturing locus in 2025: With a comparison of 'Made-in-China 2025' and 'Industry 4.0,'" *Technol. Forecast. Soc. Change*, vol. 135, no. February 2017, pp. 66–74, 2018.

[12]  L. Da Xu, E. L. Xu, and L. Li, "Industry 4.0: State of the art and future trends," *Int. J. Prod. Res.*, vol. 56, no. 8, pp. 2941–2962, 2018.

[13]  E. Armengaud, C. Sams, G. Von Falck, G. List, C. Kreiner, and A. Riel, "Industry 4.0 as digitalization over the entire product lifecycle: Opportunities in the automotive domain," *Syst. Softw. Serv. Process Improv. EuroSPI 2017. Commun. Comput. Inf. Sci.*, vol. 748, 2017.

[14]  S. Vaidya, P. Ambad, and S. Bhosle, "Industry 4.0 – A Glimpse," *Procedia Manuf.*, vol. 20, pp. 233–238, Jan. 2018.

[15]  A. Gilchrist, *Introducing Industry 4.0*. Springer, 2016.

[16]  B. Vladimirovich Sokolov, D. Ivanov, and B. Sokolov, "Integrated scheduling of material flows and information services in Industry 4.0 supply networks," *IFAC-PapersOnLine*, vol. 48, no. 3, pp. 1533–1538, 2017.

[17]  S. S. Kamble, A. Gunasekaran, and S. A. Gawankar, "Sustainable Industry 4.0 framework: A systematic literature review identifying the current trends and future perspectives," *Process Saf. Environ. Prot.*, vol. 117, pp. 408–425, Jul. 2018.

[18]  F. Shrouf, J. Ordieres, and G. Miragliotta, "Smart factories in Industry 4.0: A review of the concept and of energy management approached in production based on the Internet of

Things paradigm," in *IEEE International Conference on Industrial Engineering and Engineering Management*, 2014, pp. 697–701.

[19] R. Schmidt, M. Möhring, R.-C. Härting, C. Reichstein, P. Neumaier, and P. Jozinović, "Industry 4.0-Potentials for creating smart products: Empirical research results," *18th Int. Conf. Bus. Inf. Syst.*, vol. 208, pp. 16–27, 2015.

[20] A. C. Pereira and F. Romero, "A review of the meanings and the implications of the Industry 4.0 concept," *Procedia Manuf.*, vol. 13, pp. 1206–1214, Jan. 2017.

[21] H. S. Kang *et al.*, "Smart manufacturing: Past research, present findings, and future directions," *Int. J. Precis. Eng. Manuf. Technol.*, vol. 3, no. 1, pp. 111–128, 2016.

[22] M. Yli-Ojanperä, S. Sierla, N. Papakonstantinou, and V. Vyatkin, "Adapting an agile manufacturing concept to the reference architecture model industry 4.0: A survey and case study," *J. Ind. Inf. Integr.*, vol. 15, no. November 2018, pp. 147–160, 2019.

[23] D. Gürdür and F. Asplund, "A systematic review to merge discourses: Interoperability, integration and cyber-physical systems," *J. Ind. Inf. Integr.*, vol. 9, no. October, pp. 14–23, 2018.

[24] H. Chen, "Theoretical Foundations for Cyber-Physical Systems: A Literature Review," *J. Ind. Integr. Manag.*, vol. 02, no. 03, p. 1750013, Sep. 2017.

[25] C. Zhang, X. Xu, and H. Chen, "Theoretical foundations and applications of cyber-physical systems: a literature review," *Libr. Hi Tech*, vol. 38, no. 1, pp. 95–104, Dec. 2019.

[26] Y. Lu, "Cyber Physical System (CPS)-Based Industry 4.0: A Survey," *J. Ind. Integr. Manag.*, vol. 02, no. 03, p. 1750014, Sep. 2017.

[27] Y. Kayikci, "Sustainability impact of digitization in logistics," *Procedia Manuf.*, vol. 21,

pp. 782–789, Jan. 2018.

[28]  L. Barreto, A. Amaral, and T. Pereira, "Industry 4.0 implications in logistics: An overview," *Procedia Manuf.*, vol. 13, pp. 1245–1252, Jan. 2017.

[29]  C. Santos, A. Mehrsai, A. C. Barros, M. Araújo, and E. Ares, "Towards Industry 4.0: an overview of European strategic roadmaps," *Procedia Manuf.*, vol. 13, pp. 972–979, 2017.

[30]  S. F. Ghannadpour, S. Noori, and R. Tavakkoli-Moghaddam, "A multi-objective vehicle routing and scheduling problem with uncertainty in customers' request and priority," *J. Comb. Optim.*, vol. 28, no. 2, pp. 414–446, 2014.

[31]  V. Pillac, M. Gendreau, C. Guéret, A. Medaglia, C. Gú, and A. L. Medaglia, "A review of dynamic vehicle routing problems," *Eur. J. Oper. Res.*, vol. 225, no. 1, pp. 1–11, 2013.

[32]  G. B. Dantzig and J. H. Ramser, "The truck dispatching problem," *Manage. Sci.*, vol. 6, no. 1, pp. 80–91, Oct. 1959.

[33]  H. N. Psaraftis, M. Wen, and C. A. Kontovas, "Dynamic vehicle routing problems: Three decades and counting," *Networks An Int. J.*, vol. 67, no. 1, pp. 3–31, 2016.

[34]  N. H. M. Wilson and N. J. Colvin, *Computer control of the Rochester Dial-A-Ride System*. Cambridge: Massachusetts Institute of Technology, Center for Transportation Studies, 1977.

[35]  M. Gath, "Dispatching problems in transport logistics," in *Optimizing Transport Logistics Processes with Multiagent Planning and Control*, Wiesbaden: Springer Fachmedien Wiesbaden, 2016, pp. 15–34.

[36]  H. N. Psaraftis, "A dynamic programming solution to the single vehicle many-to-many immediate request dial-a-ride problem," *Transp. Sci.*, vol. 14, no. 2, pp. 130–154, 1980.

[37]  G. Desaulniers, "Branch-and-price-and-cut for the split-delivery vehicle routing problem

with time windows," *Oper. Res.*, vol. 58, no. 1, pp. 179–192, Feb. 2010.

[38]  P. Toth and D. Vigo, "Models, relaxations and exact approaches for the capacitated vehicle routing problem," *Discret. Appl. Math.*, vol. 123, no. 1–3, pp. 487–512, Nov. 2002.

[39]  S. Ropke and J.-F. Cordeau, "Branch and cut and price for the pickup and delivery problem with time windows," *Transp. Sci.*, vol. 43, no. 3, pp. 267–286, Aug. 2009.

[40]  M. Dror and P. Trudeau, "Stochastic vehicle routing with modified savings algorithm," *Eur. J. Oper. Res.*, vol. 23, no. 2, pp. 228–235, 1986.

[41]  G. W. Nurcahyo, R. A. Alias, and S. M. Shamsuddin, "Sweep algorithm in vehicle routing problem for public transport," *J. Antarabangsa*, vol. 2, pp. 51–64, 2002.

[42]  M. R. Khouadjia, B. Sarasola, E. Alba, L. Jourdan, and E.-G. Talbi, "A comparative study between dynamic adapted PSO and VNS for the vehicle routing problem with dynamic requests," *Appl. Soft Comput.*, vol. 12, pp. 1426–1439, 2011.

[43]  Y. Kergosien, C. Lenté, D. Piton, and J.-C. Billaut, "A tabu search heuristic for the dynamic transportation of patients between care units," *Eur. J. Oper. Res.*, vol. 214, no. 2, pp. 442–452, Oct. 2011.

[44]  M. Gendreau, F. Guertin, J.-Y. Potvin, and R. Séguin, "Neighborhood search heuristics for a dynamic vehicle dispatching problem with pick-ups and deliveries," *Transp. Res. Part C*, vol. 14, pp. 157–174, 2006.

[45]  M. Elhassania, B. Jaouad, and E. A. Ahmed, "Solving the dynamic vehicle routing problem using genetic algorithms," in *2014 International Conference on Logistics Operations Management*, 2014, pp. 62–69.

[46]  F. T. Hanshar and B. M. Ombuki-Berman, "Dynamic vehicle routing using genetic

algorithms," *Appl. Intell.*, vol. 27, no. 1, pp. 89–99, Jun. 2007.

[47]  A. E. Rizzoli, F. Oliverio, R. Montemanni, and L. M. Gambardella, "Ant Colony

Optimisation for vehicle routing problems: From theory to applications," *Swarm Intell.*,

vol. 1, no. 2, pp. 135–151, 2007.

[48]  H. Xu, P. Pu, and F. Duan, "Dynamic vehicle routing problems with enhanced ant colony

optimization," *Discret. Dyn. Nat. Soc.*, vol. 2018, pp. 1–13, Feb. 2018.

[49]  S. Chen, R. Chen, J. Gao, S. Chen, R. Chen, and J. Gao, "A Monarch Butterfly

Optimization for the Dynamic Vehicle Routing Problem," *Algorithms*, vol. 10, no. 3, p.

107, Sep. 2017.

[50]  E. Osaba, X.-S. Yang, · Fernando Diaz, E. Onieva, A. D. Masegosa, and · Asier Perallos,

"A discrete firefly algorithm to solve a rich vehicle routing problem modelling a newspaper distribution system with recycling policy," vol. 21, pp. 5295–5308, 2017.

[51]  M. Elhassania, B. Jaouad, and E. Alaoui Ahmed, "A new hybrid algorithm to solve the

vehicle routing problem in the dynamic environment," *Int. J. Soft Comput.*, vol. 8, no. 5,

pp. 327–334, 2013.

[52]  C. Novoa and R. Storer, "An approximate dynamic programming approach for the vehicle

routing problem with stochastic demands," *Eur. J. Oper. Res.*, vol. 196, no. 2, pp. 509–

515, Jul. 2009.

[53]  T. C. Du, E. Y. Li, and D. Chou, "Dynamic vehicle routing for online B2C delivery," *Int.

J. Manag. Sci.*, vol. Omega 33, pp. 33–45, 2005.

[54]  A. M. F. M. AbdAllah, D. L. Essam, and R. A. Sarker, "On solving periodic re-

optimization dynamic vehicle routing problems," *Appl. Soft Comput.*, vol. 55, pp. 1–12,

Jun. 2017.